\documentclass[amssymb,11pt]{amsart}
\usepackage{amscd,amssymb,euscript,amsthm}
\theoremstyle{plain}
\newtheorem{thm}{Theorem}[section] 
\newtheorem{cor}[thm]{Corollary}
\newtheorem{prop}[thm]{Proposition}

\theoremstyle{definition}
\newtheorem{defn}[thm]{Definition}
\theoremstyle{remark}
\newtheorem{rem}[thm]{Remark}

\numberwithin{equation}{section}
\newcommand{\id}{\operatorname{id}}

\newcommand{\ball}{\operatorname{ball}}

\def\<{\left<}
\def\>{\right>}
\def\cstar{$C^*$-algebra}
\begin{document}
\title{Lifting endomorphisms to automorphisms}
\author{William Arveson and Dennis Courtney}
%
%
\address{Department of Mathematics,
University of California, Berkeley, CA 94720}
\begin{abstract}
Normal endomorphisms of von Neumann algebras need not be extendable 
to automorphisms of a larger von Neumann algebra, but they always 
have asymptotic lifts.  We describe the structure of endomorphisms 
and their asymptotic lifts in some detail, and apply those results 
to complete the identification of asymptotic lifts of unital completely 
positive linear maps on von Neumann algebras in terms of their minimal dilations 
to endomorphisms.  
\end{abstract}
\email{arveson@math.berkeley.edu, djc@math.berkeley.edu}
\maketitle

\section{Introduction}\label{S:in}

We work in the category whose objects are pairs 
$(M,\alpha)$ consisting of a normal unit-preserving 
$*$-endomorphism $\alpha:M\to M$ of a von Neumann algebra $M$,  and whose 
maps are equivariant normal  $*$-homomorphisms that map unit to unit.  
The isomorphisms of 
this category are conjugacies, in which $\alpha_1:M_1\to M_1$ is said to be {\em conjugate} 
to $\alpha_2:M_2\to M_2$ if there is a $*$-isomorphism $\theta: M_1\to M_2$ 
satisfying $\theta\circ\alpha_1=\alpha_2\circ\theta$.  

Consider  the 
problem of extending an endomorphism $\alpha:M\to M$ to a $*$-automorphism of a larger von Neumann 
algebra, assuming that the necessary condition $\ker\alpha=\{0\}$ is satisfied.  
In that case $\alpha$ is an 
isometric $*$-endomorphism of $M$, and a straightforward construction produces 
a unital \cstar\ $N\supseteq M$ and a $*$-automorphism $\beta$ of $N$ 
that restricts to $\alpha$ on $M$.  This extension of $\alpha$ 
to an automorphism of a larger \cstar\ is unique up to natural isomorphism 
provided one assumes that it is {\em minimal} in the sense that 
$M\cup\beta^{-1}(M)\cup\beta^{-2}(M)\cup\cdots$ is  norm-dense in $N$.  

This procedure is effective for extending 
endomorphisms of \cstar s.  But it is poorly suited to this category since there is no 
{\em natural} way\footnote{It is always {\em possible} to 
carry out such a completion, 
but that construction does not give rise to 
a functor from injective endomorphisms to W$^*$-dynamical systems.  
\cite{arvKish} addresses the existence issue for $E_0$-semigroups 
acting on von Neumann algebras.} of completing the \cstar\ $N$ to a von Neumann 
algebra 
so as to obtain a $W^*$-dynamical system that 
extends $\alpha$  except in special circumstances - 
the most natural circumstance being that in which 
$\alpha$ preserves a faithful normal state of $M$.    
More serious problems arise when $\ker\alpha\neq\{0\}$, since 
in that case even extensions to $C^*$-algebraic automorphisms cannot exist.  
The proper way to associate a $W^*$-dynamical system to an endomorphism  
involves the notion of 
{\em lifting}, a concept introduced in \cite{arvLift} 
for the broader category of unital 
completely positive maps, and which will be described momentarily.  
While liftings are typically much ``smaller" than the extensions of 
isometric endomorphisms described 
above, they always exist within the 
category and they enjoy good functorial properties.  

We were led to these issues by a problem involving the broader 
category of normal unit-preserving completely positive linear maps 
$\phi:M\to M$ of von Neumann algebras $M$ (UCP maps).  
It was shown in \cite{arvLift} that 
every UCP map $\phi:M\to M$ has an asymptotic lift, 
which is unique up to natural isomorphism.  Naturally, one wants to identify 
the asymptotic lift of $\phi$ in concrete terms.  In \cite{arvLift},  
the asymptotic lift of $\phi$ was identified 
as the tail flow of the minimal dilation of $\phi$ 
in ``most" cases - namely those cases in which 
the tail flow of the dilated endomorphism has trivial kernel.  But in general, the 
minimal dilation of $\phi$ 
to an endomorphism can have a nontrivial kernel, and 
the identification 
problem was left open in those cases.  

The purpose of this note is to identify asymptotic lifts of UCP maps 
on von Neumann algebras in general.  
That is accomplished by first giving a description of lifts of endomorphisms, 
in the course of which we obtain a basic result on the structure 
of surjective endomorphisms of von Neumann algebras 
that appears to have been overlooked (Theorem \ref{seThm1}).  We apply these results 
to identify the asymptotic lift of 
an arbitrary UCP map in terms of its minimal dilation to an endomorphism 
of a larger von Neumann algebra (Theorem \ref{apThm1}), thereby completing 
Theorem 7.1 of \cite{arvLift}.  

In related work \cite{arvStor}, the notion of asymptotic lift 
was generalized to normal positive linear maps 
acting on von Neumann algebras (also see \cite{storPos}).  
It is significant that 
since there is no dilation theory for positive linear maps that 
are not completely positive,  the identification 
problem becomes a significant issue in such cases and has been only partially solved.  
Further discussion can be found in \cite{arvStor}.

\section{Lifting endomorphisms}\label{S:le}

Throughout this section, $\alpha:M\to M$ will denote an endomorphism 
acting on a von Neumann algebra $M$.  
By a {\em $W^*$-dynamical system} 
we mean a pair $(N,\beta)$, where $\beta$ is a $*$-automorphism 
of a von Neumann algebra $N$.

\begin{defn}
A {\em lifting} of $\alpha:M\to M$ is a triple $(N, \beta,E)$ where 
$(N,\beta)$ is a $W^*$-dynamical system and $E:N\to M$ is a unit-preserving 
normal $*$-homomorphism satisfying $E\circ\beta=\alpha\circ E$.  
\end{defn}

Note first that for every lifting $(N,\beta,E)$ of $\alpha$, we have 
\begin{equation}
E(N)\subseteq M\cap\alpha(M)\cap\alpha^2(M)\cap\cdots.  
\end{equation}
Indeed, every element $y=E(x)$ in the range of $E$ can be written in 
the form $y=\alpha^n(E(\beta^{-n}(x))\in\alpha^n(M)$ for every $n=0,1,2,\dots$, 
from which the assertion is evident.  

\begin{rem}[Nondegeneracy]\label{leRem1}
A lifting $(N,\beta,E)$ of $\alpha$ is said to be {\em nondegenerate} 
if for every $x\in N$, 
$$
E(\beta^n(x))=0, \quad n\in \mathbb Z\quad\implies x=0.  
$$
In general, the set 
$$
K=\{x\in N: E(\beta^n(x))=0,\quad\forall\ n\in\mathbb Z\}
$$
is a weak$^*$-closed two-sided ideal in $N$ satisfying $\beta(K)=K$.  Hence 
there is a $\beta$-fixed central projection $c\in N$ such that $K=cN$.  It 
follows that $N$ decomposes into a sum  
$N=K\oplus N_0$, where  $(N_0,\beta\restriction_{N_0}, E\restriction_{N_0})$ 
is a nondegenerate lift of $\alpha$ and $E(\beta^n(K))=\{0\}$ for every $n\in \mathbb Z$.  
In particular, {\em every lift $(N,\beta,E)$ of $\alpha$ can be 
reduced to a nondegenerate lift without affecting the range of the homomorphism 
$E$.  }
\end{rem}

Let $\alpha:M\to M$ be an endomorphism of a von Neumann algebra.  The 
sequence of von Neumann algebras $M, \alpha(M), \alpha^2(M),\dots$ decreases as $n$ increases, 
and their intersection 
$$
M_\infty=\bigcap_{n=1}^\infty\alpha^n(M)
$$
is called the {\em tail algebra} of $\alpha$. The restriction of $\alpha$ 
to the tail algebra is a {\em surjective} endomorphism; it is an automorphism 
iff $\ker\alpha\cap M_\infty=\{0\}$.

\begin{prop}\label{leProp1}  Let $\alpha:M\to M$ be an endomorphism.  
For every lifting $(N,\beta,E)$ of $\alpha$, the following are equivalent. 
\begin{enumerate}
\item[(i)] For every normal linear functional $\rho\in M_*$, one has 
\begin{equation}\label{leEq2}
\lim_{n\to\infty}\|\rho\circ\alpha^n\|=\|\rho\circ E\|.  
\end{equation}
\item[(ii)] $E(N)=M_\infty$.  
\end{enumerate}
\end{prop}

\begin{proof}
Since $E$ is a $*$-homomorphism of von Neumann algebras, it maps the unit ball of $N$ onto the unit 
ball of its range.  Hence (ii) is equivalent to 
$$
E(\ball N)=\ball M_\infty=\bigcap_{n=1}^\infty\alpha^n(\ball M).
$$  
The equivalence (i)$\iff$(ii) now follows from the more general assertion of 
Lemma 3.6 of \cite{arvLift}.
\end{proof}

\begin{defn}\label{leDef2}
An {\em asymptotic lift} of an endomorphism $\alpha:M\to M$ is a nondegenerate lifting $(N,\beta,E)$ 
satisfying 
the conditions of Proposition \ref{leProp1}.  
\end{defn}

\begin{rem}
[Relation to asymptotic lifts of UCP maps]  In \cite{arvLift}, the term {\em asymptotic lift} refers to a related 
concept that was introduced for the broader category of UCP maps on dual operator systems.  
It is significant that an asymptotic lift in the sense of Definition \ref{leDef2} is 
also an asymptotic lift in the broader sense of Definition 3.1 of \cite{arvLift}.  

To prove that assertion, it suffices 
to show that if a lifting $(N,\beta,E)$ of an endomorphism $\alpha:M\to M$ satisfies 
the equivalent properties (i) and  (ii) of Proposition \ref{leProp1}, then those properties 
persist throughout the matrix hierarchy over $M$.  Indeed, for each $n=1,2,\dots$, 
the lifting $(N,\beta,E)$ of $\alpha:M\to M$ induces a natural lift $(M_n\otimes N,\id_n\otimes\beta, \id_n\otimes E)$ of 
the endomorphism $\id_n\otimes\alpha:M_n\otimes M\to M_n\otimes M$, 
and by examining matrix entries 
one finds that property (ii) persists at level $n$.  Hence property (i) holds as well for 
every $n=1,2,\dots$, and $(N,\beta,E)$ satisfies Definition 3.1 of \cite{arvLift}.  
\end{rem}

Two liftings $(N_k,\beta_k,E_k)$, $k=1,2$,  of an endomorphism $\alpha:M\to M$ 
are said to be {\em isomorphic } if there is an isomorphism of von Neumann algebras 
$\gamma: N_1\to N_2$ satisfying $\gamma\circ\beta_1=\beta_2\circ\gamma$ and 
$E_2\circ\gamma=E_1$.

\begin{thm}\label{leThm1}
Every endomorphism $\alpha:M\to M$  of a von Neumann algebra has an asymptotic  
lifting that is unique up to isomorphism.   
\end{thm}

\begin{proof}
The proof of Theorem 3.2 of \cite{arvLift} explicitly constructs an asymptotic lift of a UCP 
map of that category in terms of the space of inverse sequences of that map.  Since in  
the present context the 
map is an endomorphism $\alpha:M\to M$, one sees by inspection that the constructed asymptotic 
lift $(N,\beta, E)$ has the following properties:  the space 
$N$ of inverse sequences is closed under multiplication, $\beta$ 
is a $*$-automorphism of that von Neumann algebra, and $E:N\to M$ is a normal $*$-homomorphism.  Hence 
$(N,\beta,E)$ is an asymptotic lift in the sense of Definition \ref{leDef2}.  
The proof of uniqueness involves similar observations.   
\end{proof}

\section{Structure of surjective endomorphisms}\label{S:se}

In this section we prove 
that in general, a surjective endomorphism of a von Neumann algebra 
admits a natural decomposition into the direct sum of a $W^*$-dynamical system and 
an endomorphism of a particularly simple kind, called a {\em backward shift},  that 
depends only on  $\ker\alpha$.  That allows us to identify 
asymptotic lifts of endomorphisms in very concrete terms.

Let $K$ be a von Neumann algebra and consider the von Neumann algebra 
$\ell^\infty(\mathbb N,K)$ of all singly-infinite bounded sequences 
$$
x=(x_1, x_2,\dots),\qquad x_k\in K.  
$$
Define an endomorphism $\sigma_+$ of $\ell^\infty(\mathbb N,K)$ as follows
\begin{equation}\label{seEq1}
\sigma_+(x_1,x_2,\dots)=(x_2,x_3,\dots),\qquad (x_1,x_2,\dots)\in \ell^\infty(\mathbb N,K).  
\end{equation}
Obviously, $\sigma_+$ is a normal surjective unit-preserving endomorphism, and 
$$
\ker\sigma_+=K\oplus 0\oplus 0\oplus \cdots \cong K.   
$$
Such an endomorphism $\sigma_+$ is called the {\em backward shift} based on $K$.

We can modify the backward shift $\sigma_+$ based on $K$ in a nontrivial 
way by choosing an automorphism $\beta$ of 
another von Neumann algebra $P$ and letting $\beta\oplus\sigma_+$ be the 
endomorphism of $P\oplus \ell^\infty(\mathbb N,K)$ defined by 
\begin{equation}\label{seEq2}
\beta\oplus\sigma_+: x\oplus y\in P\oplus \ell^\infty(\mathbb N,K)\mapsto \beta(x)\oplus\sigma_+(y).  
\end{equation}
This is a surjective endomorphism whose kernel is isomorphic 
to $K$, but it has a summand $P$ on which it restricts to 
an automorphism.  

\begin{thm}\label{seThm1}
Every normal surjective endomorphism $\alpha$ of a von Neumann algebra is conjugate 
to one of the form (\ref{seEq2}), where $\sigma_+$ is the backward shift 
based on $K=\ker \alpha$,  and where $\beta$ is the $*$-automorphism  
defined by restricting $\alpha$ to 
$P=(\mathbf 1-c) M$, $c$ being the $\alpha$-fixed central projection  
\begin{equation}\label{seEq3}
c=\lim_{n\to\infty}c_n,  
\end{equation}
where $c_1\leq c_2\leq\cdots$ is the sequence of central projections 
$\ker\alpha^n=c_n M$.     
\end{thm}

\begin{proof}
Since  $\ker\alpha^n$ is a weak$^*$-closed ideal in $M$ for every $n=1,2,\dots$, it 
has the form $c_nM$ where 
$c_n$ is a central projection; and since $\ker\alpha^n\subseteq \ker\alpha^{n+1}$, 
it follows that $c_n\leq c_{n+1}$.

We claim: $\alpha(c_{n+1})=c_n$, for every $n=1,2,\dots$.  
Since $\ker\alpha^k=c_kM$ and $\alpha(M)=M$, 
this is equivalent to the assertion $\alpha(\ker\alpha^{n+1})=\ker\alpha^n$.
Obviously,  $x\in \ker\alpha^{n+1}\implies\alpha(x)\in \ker\alpha^n$.  
 For the opposite inclusion, choose $y\in\ker\alpha^n$.  
Since $M=\alpha(M)$, we can find $x\in M$ such that $y=\alpha(x)$.  This $x$ must satisfy 
$
\alpha^{n+1}(x)=\alpha^n(\alpha(x))=\alpha^n(y)=0, 
$ 
hence $x\in \ker\alpha^{n+1}$, and therefore $y=\alpha(x)\in\alpha(\ker\alpha^{n+1})$.  
These formulas  
$\alpha(c_{n+1})=c_{n}$ clearly imply that the limit $c$ of (\ref{seEq3}) is
an $\alpha$-fixed central projection.

Let $P=(\mathbf 1-c) M$.  Note that $\alpha(P)=\alpha(\mathbf 1-c)\alpha(M)=(\mathbf 1-c)M=P$.  
Since the kernel of $\alpha$ is $c_1 M\subseteq c M$ we have 
$P\cap\ker\alpha=\{0\}$, hence $\alpha$ restricts to a $*$-automorphism of $P$.  

Turning now to the summand $c M$, we claim that for every $n\geq 2$, 
$\alpha^{n-1}$ restricts to an isomorphism of von Neumann algebras 
$$
\alpha^{n-1}:(c_n-c_{n-1})M\cong \ker \alpha=c_1 M.  
$$  
Indeed, the restriction of $\alpha^{n-1}$ to $(c_n-c_{n-1})M$ is injective 
because $\ker\alpha^{n-1}=c_{n-1}M$ intersects trivially with the 
algebra on the left.  It is surjective because after iterating the 
formulas $\alpha(c_{k+1})=c_k$ we find that 
\begin{align*}
\alpha^{n-1}((c_n-c_{n-1})M)&=\alpha^{n-1}(c_nM)=\alpha^{n-1}(c_n)M
=\alpha^{n-2}(c_{n-1})M
\\
&=\cdots=\alpha(c_2)M=c_1M.
\end{align*}

Now consider the von Neumann algebra $\ell^\infty(\mathbb N,\ker\alpha)$, 
the algebra of all uniformly bounded sequences $y=(y_1,y_2,\dots)$ with 
$y_k\in \ker\alpha$ for $k\geq 1$.   Every element $x\in c M$ admits 
a unique decomposition into a bounded sequence of mutually orthogonal central slices 
$
x=x_1+x_2+x_3+\cdots, 
$
where $x_1=c_1x$ and $x_k=(c_k-c_{k-1})x$ for $k\geq 2$.  
Moreover, the preceding paragraph implies that $\alpha^{k-1}(x_k)\in\ker\alpha$ 
for every $k\geq 2$.  Thus we can define a normal homomorphism of von Neumann 
algebras $\theta: cM\to \ell^\infty(\mathbb N,\ker\alpha)$ by 
$$
\theta(x)=(x_1,\alpha(x_2),\alpha^2(x_3),\dots),\qquad x\in c M.  
$$
We have also seen that for each $k\geq 2$, $\alpha^{k-1}$ restricts to 
an isomorphism from $(c_k-c_{k-1})M$ to $\ker \alpha$; and since 
$cM=c_1M\oplus (c_2-c_1)M\oplus (c_3-c_2)M\oplus\cdots$, 
it follows that $\theta$ is an isomorphism 
of von Neumann algebras.  

One can now directly verify that $\theta\circ\alpha=\sigma_+\circ\theta$, where 
$\sigma_+$ denotes the backward shift on $\ell^\infty(\ker\alpha)$.  We conclude that  
$\theta$ implements a conjugacy of the restriction of $\alpha$ to $c M$ 
and the backward shift based on $\ker\alpha$.  
\end{proof}

The $W^*$-dynamical system $(P,\alpha\restriction_P)$ is called the 
{\em automorphic summand} of $\alpha$.  
It is clear from the preceding proof that 
{\em two surjective endomorphisms are conjugate iff 
their automorphic summands are conjugate $W^*$-dynamical systems 
and their kernels are isomorphic von Neumann algebras.}

Theorem \ref{seThm1} leads to the following description of the asymptotic lifts 
of arbitrary endomorpisms of von Neumann algebras.  

\begin{cor}\label{seCor1}  Let $\alpha:M\to M$ be an endomorphism with tail algebra 
$M_\infty=\cap_{n\geq 1}\alpha^n(M)$.  Let $K=\ker \alpha\cap M_\infty$, 
let $\sigma_+$ be the backward shift acting on $\ell^\infty(\mathbb N,K)$,  
and let $(P,\beta)$ be the automorphic summand of $\alpha\restriction_{M_\infty}$.  

By Theorem \ref{seThm1}, there is an isomorphism of von Neumann algebras 
$$
\theta: P\oplus\ell^\infty(\mathbb N,K)\to M_\infty
$$ 
that satisfies $\theta\circ (\beta\oplus\sigma_+)=\alpha\circ\theta$.  Let 
$\sigma$ be the bilateral shift acting on the von Neumann algebra 
$\ell^\infty(\mathbb Z,K)$ by way of $\sigma(x_n)=(x_{n+1})$, 
and define a homomorphism $E: P\oplus \ell^\infty(\mathbb Z,K)\to M_\infty$ by 
$$
E(p\oplus (x_n))=\theta(p\oplus (x_1,x_2,\dots)), \qquad p\in P,\quad (x_n)\in\ell^\infty(\mathbb Z,K). 
$$
Then $(P\oplus\ell^\infty(\mathbb Z,K),\beta\oplus \sigma,E)$ is the asymptotic 
lift of $\alpha:M\to M$.  
\end{cor}

\begin{proof}
It is obvious that $(P\oplus\ell^\infty(\mathbb Z,K), \beta\oplus\sigma)$ is a W$^*$-dynamical system and 
that $E$ is a homomorphism of von Neumann algebras with range 
$$
E(P\oplus\ell^\infty(\mathbb Z,K))=\theta(P\oplus \ell^\infty(\mathbb N,K))=M_\infty.  
$$
Moreover, 
\begin{align*}
E\circ(\beta\oplus\sigma)(p\oplus(x_n))&=E(\beta(p)\oplus (x_{n+1}))=\theta(\beta(p)\oplus(x_2,x_3,\dots))
\\
&
=
\theta(\beta(p)\oplus\sigma_+(x_1,x_2,\dots))=\alpha\circ\theta(p\oplus(x_1,x_2,\dots))
\\
&
=\alpha\circ E(p\oplus (x_n)), 
\end{align*}
hence $E\circ(\beta\oplus\sigma)=\alpha\circ E$.  We conclude that 
$(P\oplus\ell^\infty(\mathbb Z,K), \beta\oplus\sigma,E)$ is a lifting of $\alpha$ that 
satisfies condition (ii) of Proposition \ref{leProp1}; and it remains only to show that this 
lifting is nondegenerate.  But if $p\in P$ and $(x_n)\in\ell^\infty(\mathbb Z,K)$ are such that 
$E((\beta\circ\sigma)^k(p\oplus(x_n))=\theta(\beta^k(p)\oplus (x_{k+1},x_{k+2},\dots))=0$ 
for every $k\in \mathbb Z$, then $\beta^k(p)=0$ and 
$x_{k+1}=0$ for every $k\in\mathbb Z$. The desired formula $p\oplus(x_n)=0$ follows.  
\end{proof}

\section{Application to UCP maps on von Neumann algebras}\label{S:ap}

Let $\phi:M\to M$ be a UCP map acting on a von Neumann algebra $M$.  
In this section we identify the asymptotic lift of $\phi$ 
in terms of its minimal dilation to an endomorphism of a larger 
von Neumann algebra.  This solves the identification problem in general by
strengthening Theorem 7.1 of \cite{arvLift} that was restricted to the case 
in which the minimal dilation has trivial kernel.  Indeed, the following result 
applies to dilations of $\phi$ that are not necessarily minimal (see Chapter 
8 of \cite{arvMono}).  

\begin{thm}\label{apThm1}
Let $\alpha:N\to N$ be an endomorphism of a von Neumann algebra 
and let $p\in N$ be a projection that satisfies $\alpha(p)\geq p$ and 
$\alpha^n(P)\uparrow\mathbf 1$ as $n\uparrow\infty$.  Let $M=pNp$ 
and let $\phi:M\to M$ be the UCP map defined by 
$$
\phi(x)=p\alpha(x)p,\qquad x\in M=pNp.  
$$
Let $(\tilde N,\tilde\alpha, E)$ be the asymptotic lift of 
$\alpha$ described in Corollary \ref{seCor1}.  Then the asymptotic lift 
of $\phi$ is  $(\tilde N,\tilde\alpha,\tilde E)$ where 
$\tilde E: \tilde N\to M$ is the UCP map 
\begin{equation}\label{apEq1}
\tilde E(x)=p E(x)p,\qquad x\in \tilde N.   
\end{equation}
\end{thm}

\begin{proof} Obviously $\tilde E:\tilde N\to M$ is a UCP map and 
we claim $\phi\circ \tilde E=\tilde E\circ \tilde\alpha$.  Indeed, 
we can use $p\alpha(p)=\alpha(p)p=p$ and $\alpha\circ E=E\circ\tilde \alpha$ to write 
$$
\phi(\tilde E(x))=p\alpha(\tilde E(x))p=p\alpha(p)\alpha(E(x))\alpha(p)p=
pE(\tilde\alpha(x))p=\tilde E\circ\tilde\alpha(x).  
$$
Hence $(\tilde N,\tilde\alpha,\tilde E)$ is a lifting of $\phi$.  To see that it 
is nondegenerate, choose $x\in \tilde N$ such that $\tilde E(\tilde\alpha^k(x))=0$, $k\in\mathbb Z$.  
Then for $n\geq 1$ we can apply $\alpha^n$ to $\tilde E(\tilde\alpha^{-n}(x))=0$ and 
use $\alpha\circ E=E\circ\tilde\alpha$ to obtain 
$$
0=\alpha^n(\tilde E(\tilde\alpha^{-n}(x)))=\alpha^n(p)\alpha^n(E(\tilde\alpha^{-n}(x)))\alpha^n(p)
=\alpha^n(p)E(x)\alpha^n(p).  
$$
Since $\alpha^n(p)\uparrow \mathbf 1$ as $n\uparrow\infty$, it follows that $E(x)=0$.  Replacing 
$x$ with $\tilde\alpha^k(x)$, $k\in\mathbb Z$, and using nondegeneracy of $(\tilde N,\tilde \alpha,E)$, we 
conclude that $x=0$.  

We claim that for every $\rho\in M_*$, 
\begin{equation}\label{apEq2}
\lim_{n\to\infty}\|\rho\circ\phi^n\|=\|\rho\circ\tilde E\|.  
\end{equation}
To prove (\ref{apEq2}), fix $\rho$ and define a normal functional 
$\bar\rho\in N_*$ by $\bar\rho(y)=\rho(pyp)$.  For every $x\in\tilde N$ we have 
$\rho\circ\tilde E(x)=\rho(pE(x)p)=\bar\rho\circ E(x)$, and 
as in  the proof of 
formula (7.2) of \cite{arvLift}, we obtain the following formulas for $n\geq 1$
\begin{equation}\label{apEq3}
\|\rho\circ\phi^n\|=\|\bar\rho\circ\alpha^n\|,\qquad \|\rho\circ \tilde E\|=\|\bar\rho\circ E\|.   
\end{equation}
Since $(\tilde N,\tilde\alpha,E)$ is the asymptotic lift of $\alpha:N\to N$, $\|\bar\rho\circ\alpha^n\|$ 
converges to $\|\bar\rho\circ E\|$ as $n\to\infty$, and (\ref{apEq2}) follows.  
Similarly, one can promote (\ref{apEq2}) throughout the matrix hierarchy over $M$ exactly as 
in the proof of Theorem 7.1 of \cite{arvLift} to complete the proof of Theorem \ref{apThm1}.  
\end{proof}

\bibliographystyle{alpha}

\newcommand{\noopsort}[1]{} \newcommand{\printfirst}[2]{#1}
  \newcommand{\singleletter}[1]{#1} \newcommand{\switchargs}[2]{#2#1}

\end{document}